\documentclass[10pt]{amsart}
\usepackage[cp1251]{inputenc}
\usepackage[english,russian]{babel}
\usepackage{amsmath}
\usepackage{amssymb}

\usepackage{amsfonts}
\usepackage{graphicx}

\def\udcs{???.?} 
\def\mscs{??X??} 
\def\logo{{\bf\huge S\raisebox{0.2ex}{\hspace{0.55ex}\raisebox{0.05ex}e\hspace{-1.65ex}$\bigcirc$}MR}}

\def\semrtop
     {
  \vbox{
     \noindent\logo
     \hspace{80mm}\raisebox{1ex}{ISSN 1813-3304 }

     \vspace{5mm}

     \begin{center}
     {\huge СИБИРСКИЕ \ ЭЛЕКТРОННЫЕ} \\[2mm]  
     {\huge МАТЕМАТИЧЕСКИЕ ИЗВЕСТИЯ} \\[2mm]
     {\large Siberian Electronic Mathematical Reports} \\[1mm]
     {\LARGE\tt{http://semr.math.nsc.ru}}\\[0.5mm]
     \end{center}
     \vspace{-3mm}
     \noindent
     \begin{tabular}{c}
     \hphantom{aaaaaaaaaaaaaaaaaaaaaaaaaaaaaaaaaaaaaaaaaaaaaaaaaaaaaaaaaaaaaaaaaaaaaa} \\
     \hline\hline
     \end{tabular}

     \vspace{1mm}
     {\flushleft\it Том 16, стр. 144--144 (2019) \hspace{65mm}{\rm\small УДК \udcs}} 
     \newline
         {\rm\small DOI~10.33048/semi.2019.16.xxx}\hphantom{aaaaaaaaaaaaaaaaaaaaaaaaaaaaaaaaaaaa}{\rm\small MSC\ \ \mscs }
  }
}

\numberwithin{equation}{section}
\newtheorem{thm}{Theorem}[section]

\newtheorem{lem}[thm]{Lemma}

\newtheorem{re}{Remark}[section]

\newenvironment{pf}{{\noindent \it \bf Proof:}}{{\hfill$\Box$}\\}

\begin{document}
\renewcommand{\refname}{References}
\renewcommand{\proofname}{Proof.}
\renewcommand{\figurename}{Fig.}

\title{A regularity criterion to the 3D Boussinesq equations}

\author{{A.M. Alghamdi}}%
\address{Ahmad Mohammad Alghamdi
\newline\hphantom{iii} Department of Mathematical Science , Faculty of Applied Science,
\newline\hphantom{iii} Umm Alqura University, P.O.B. 14035, Makkah 21955, Saudi Arabia}%
\email{amghamdi@uqu.edu.sa; }%

\author{{I. Ben Omrane}}%
\address{Ines Ben Omrane
\newline\hphantom{iii} Department of Mathematics, Faculty of Science, Imam Mohammad Ibn Saud,
\newline\hphantom{iii} Islamic University (IMSIU), P. O. Box 90950, Riyadh 11623, Saudi Arabia}%
\email{:imbenomrane@imamu.edu.sa}%

\author{{S. Gala}}%
\address{Sadek Gala
\newline\hphantom{iii} Department of Mathematics, Ecole Normale Sup\'{e}rieure de Mostaganem
 University of Mostaganem,  Box 227, Mostaganem 27000, Algeria,
\newline\hphantom{iii} Dipartimento di Matematica e Informatica, Viale Andrea Doria, 6, 95125- Catania, Italy}
\email{sgala793@gmail.com}%

\author{{M.A.Ragusa}}%
\address{Maria Alessandra Ragusa
\newline\hphantom{iii} Dipartimento di Matematica e Informatica, Viale Andrea Doria, 6, 95125- Catania, Italy
\newline\hphantom{iii}  RUDN University, 6 Miklukho - Maklay St, Moscow, 117198, Russia}
\email{maragusa@dmi.unict.it}%

\thanks{\rm The researchers acknowledge the Deanship of Scientific Research at Imam
Mohammad Ibn Saud Islamic University, Saudi Arabia, for financing this
project under the grant no. (381206). Part of the work was carried out while
the third author was long-term visitor at University of Catania. The
hospitality of Catania University is graciously acknowledged. This research
is partially supported by Piano della Ricerca 2016-2018 - Linea di
intervento 2: "Metodi variazionali ed equazioni differenziali". M.A. Ragusa
wish to thank the support of "RUDN University Program 5-100". The authors
wish to express their thanks to the referees for their very careful reading
of the paper, giving valuable comments and helpful suggestions.}
\thanks{\it Received  July, 30, 2019, published  ...,  2020.}%

\date{}
\semrtop \vspace{1cm}
\maketitle {\small
\begin{quote}
\noindent{\sc Abstract. } \medskip
The paper deals with the regularity criterion for the weak solutions to the
3D Boussinesq equations in terms of the partial derivatives in Besov spaces.
It is proved that the weak solution $(u,\theta )$ becomes regular provided
that
\begin{equation*}
(\nabla _{h}u,\nabla _{h}\theta )\in L^{\frac{8}{3}}(0,T;\overset{\cdot }{B}%
_{\infty ,\infty }^{-1}(\mathbb{R}^{3})).
\end{equation*}%
Our results improve and extend the well-known results of Fang-Qian \cite{FQ}
for the Navier-Stokes equations.\medskip

 \end{quote}
}



%
\section{Introduction and main result}

This paper is devoted to the study of the Cauchy problem for the Boussinesq
equations in $\mathbb{R}^{3}\times (0,T)$:
\begin{equation}
\left\{
\begin{array}{l}
\partial _{t}u-\Delta u+u\cdot \nabla u+\nabla \pi =\theta e_{3}, \\
\partial _{t}\theta -\Delta \theta +u\cdot \nabla \theta =0, \\
\nabla \cdot u=0, \\
u\left( x,0\right) =u_{0}\left( x\right) ,\text{ \ \ \ \ }\theta \left(
x,0\right) =\theta _{0}\left( x\right) ,%
\end{array}%
\right.  \label{eq1.1}
\end{equation}%
where $u=u(x,t)$ is the velocity of the fluid, $\theta =\theta (x,t)$ is the
scalar quantity such as the concentration of a chemical substance or the
temperature variation in a gravity field, $\pi =\pi (x,t)$ is the scalar
pressure, while $u_{0}$ and $\theta _{0}$ are given initial velocity and
initial temperatture with $\nabla \cdot u_{0}=0$ in the sense of
distributions. $e_{3}=(0,0,1)^{T}$ denotes the vertical unit vector.

The Cauchy problem (\ref{eq1.1}) for the Boussinesq equation, has been
studied extensively by many authors (see, for example, \cite{A, AH, BS, can,
CN, CKN, DJZ, FO, FZ, G1, G6, G3, G7, G8, HL, IM} and references cited
therein).

When $\theta =0$, (\ref{eq1.1}) is the well-known incompressible
Navier-Stokes equations, which the global regularity is an outstanding open
problem, as well as the famous millennium prize problem. Since the global
existence of weak solutions is well-known and strong solutions are unique
and smooth in $(0,T)$, it is an interesting problem on the regularity
criterion of the weak solutions if some partial derivatives of the velocity
satisfy certain growth conditions (see, e.g. \cite{Ber, FQ, G2, G4, G5, KZ,
Ska, Z1, Z2}). One of the most significant achievements in this direction is
the celebrated Fang and Qian criterion \cite{FQ}. More precisely, they
showed that a weak solution with $H^{1}-$data is a strong solution provided
that
\begin{equation}
\nabla _{h}u\in L^{\frac{8}{3}}(0,T;\overset{\cdot }{B}_{\infty ,\infty
}^{-1}(\mathbb{R}^{3})).  \label{eq0.1}
\end{equation}%
where $\nabla _{h}=\left( \partial _{1},\partial _{2}\right) $ denotes the
horizontal gradient operator and $\overset{\cdot }{B}_{\infty ,\infty }^{-1}$
denotes the homogeneous Besov space. For details see \cite{T}.

Recall that the weak solutions satisfy the following energy inequality
\begin{equation}
\left\Vert u(t)\right\Vert _{L^{2}}^{2}+\left\Vert \theta (t)\right\Vert
_{L^{2}}^{2}+2\int_{0}^{t}\!\!(\left\Vert \nabla u(\tau )\right\Vert
_{L^{2}}^{2}+\left\Vert \nabla \theta (\tau )\right\Vert _{L^{2}}^{2})d\tau
\leq \left\Vert u_{0}\right\Vert _{L^{2}}^{2}+\left\Vert \theta
_{0}\right\Vert _{L^{2}}^{2},  \label{eq100}
\end{equation}%
for all $0\leq t\leq T$.

Motivated by the reference mentioned above, our aim of the present paper is
to improve and extend the above regularity criterion (\ref{eq0.1}) to the
Boussinesq equations (\ref{eq1.1}).

Our main result reads as follows.

\begin{thm}
\label{th1}Suppose $T>0$, $\left( u_{0},\theta _{0}\right) \in H^{1}(\mathbb{%
R}^{3})$ with $\mathrm{div}\,u_{0}=0$ in $\mathbb{R}^{3}$, in the sense of
distributions. Let $(u,\theta )$ be a weak solution of (\ref{eq1.1}) in $%
(0,T)$. Assume that%
\begin{equation}
(\nabla _{h}u,\nabla _{h}\theta )\in L^{\frac{8}{3}}(0,T;\overset{\cdot }{B}%
_{\infty ,\infty }^{-1}(\mathbb{R}^{3})),  \label{eq6}
\end{equation}%
then the weak solution $(u,\theta )$ is regular on $\mathbb{R}^{3}\times
(0,T]$.
\end{thm}

\begin{re}
In the case $\theta =0$, the above theorem reduces to the well-known Fang
and Qian result \cite{FQ} for the Navier-Stokes equations.
\end{re}

\subsection{Proof of Theorem \protect\ref{th1}}

In this section, we shall give the proof of Theorem \ref{th1}, we first need
to prove the following lemma.

\begin{lem}
\label{lem0}Let $(u,\theta )$ be a smooth solution to (\ref{eq1.1}). Then,
there exists a positive universal constant $C$ such that the following a
priori estimates hold :%
\begin{eqnarray}
&&\int_{\mathbb{R}^{3}}(u\cdot \nabla )u\cdot \Delta udx+\int_{\mathbb{R}%
^{3}}(u\cdot \nabla )\theta \cdot \Delta \theta dx  \notag \\
&\leq &C\int_{\mathbb{R}^{3}}\left\vert \nabla _{h}u\right\vert \left\vert
\nabla u\right\vert ^{2}dx+C\int_{\mathbb{R}^{3}}\left\vert \nabla
_{h}u\right\vert \left\vert \nabla \theta \right\vert ^{2}dx+C\int_{\mathbb{R%
}^{3}}\left\vert \nabla _{h}\theta \right\vert \left\vert \nabla
u\right\vert \left\vert \nabla \theta \right\vert dx.  \label{eq20.1}
\end{eqnarray}
\end{lem}

\begin{pf}
Due to the divergence-free condition $\nabla \cdot u=0$, one shows that%
\begin{eqnarray*}
\sum_{i,j,k=1}^{3}\int_{\mathbb{R}^{3}}u_{i}\partial _{i}\partial
_{k}u_{j}\partial _{k}u_{j}dx &=&\frac{1}{2}\sum_{i,j,k=1}^{3}\int_{\mathbb{R%
}^{3}}u_{i}\partial _{i}(\partial _{k}u_{j})^{2}dx \\
&=&-\frac{1}{2}\int_{\mathbb{R}^{3}}(\sum_{i=1}^{3}\partial _{i}u_{i})\left(
\sum_{j,k=1}^{3}(\partial _{k}u_{j})^{2}\right) dx=0,
\end{eqnarray*}%
\begin{equation*}
\sum_{i,j,k=1}^{3}\int_{\mathbb{R}^{3}}u_{i}\partial _{i}\partial _{k}\theta
\partial _{k}\theta dx=0,
\end{equation*}%
Hence,%
\begin{eqnarray*}
I &=&\int_{\mathbb{R}^{3}}(u\cdot \nabla )u\cdot \Delta udx\int_{\mathbb{R}%
^{3}}(u\cdot \nabla )\theta \cdot \Delta \theta dx \\
&=&-\int_{\mathbb{R}^{3}}\nabla (u\cdot \nabla )u\cdot \nabla udx-\int_{%
\mathbb{R}^{3}}\nabla (u\cdot \nabla )\theta \cdot \nabla \theta dx \\
&=&-\sum_{i,j,k=1}^{3}\int_{\mathbb{R}^{3}}\partial _{k}u_{i}\partial
_{i}u_{j}\partial _{k}u_{j}dx-\sum_{i,k=1}^{3}\int_{\mathbb{R}^{3}}\partial
_{k}u_{i}\partial _{i}\theta \partial _{k}\theta dx.
\end{eqnarray*}

In order to estimate the right hand side of $I$, we split each of the above
integrals according to the following rules :

\textbf{Case 1 }: when $1\leq k\leq 2$ or $1\leq i\leq 2$, then the integral
$I$ has at least $\nabla _{h}u$ or $\nabla _{h}\theta $ in the integrand and
can be dominated by%
\begin{equation}
I\leq C\int_{\mathbb{R}^{3}}\left\vert \nabla _{h}u\right\vert \left\vert
\nabla u\right\vert ^{2}dx+C\int_{\mathbb{R}^{3}}\left\vert \nabla
_{h}u\right\vert \left\vert \nabla \theta \right\vert ^{2}dx+C\int_{\mathbb{R%
}^{3}}\left\vert \nabla _{h}\theta \right\vert \left\vert \nabla
u\right\vert \left\vert \nabla \theta \right\vert dx  \label{eq03}
\end{equation}

\textbf{Case 2 }: when $k=i=3$, then we use the divergence-free condition to
rewrite
\begin{equation*}
\partial _{3}u_{3}=-\partial _{1}u_{1}-\partial _{2}u_{2},
\end{equation*}%
then the integral can be controlled by (\ref{eq03}). Hence the proof of
Lemma is complete.
\end{pf}

\subsection{Proof of Theorem \protect\ref{th1}}

Before going to the proof, we recall the following inequality established in
\cite{CW} :
\begin{eqnarray}
\left\Vert f\right\Vert _{L^{r}} &\leq &C\left\Vert f\right\Vert _{L^{2}}^{%
\frac{6-r}{2r}}\left\Vert \partial _{1}f\right\Vert _{L^{2}}^{\frac{r-2}{2r}%
}\left\Vert \partial _{2}f\right\Vert _{L^{2}}^{\frac{r-2}{2r}}\left\Vert
\partial _{3}f\right\Vert _{L^{2}}^{\frac{r-2}{2r}}  \notag \\
&\leq &C\left\Vert f\right\Vert _{L^{2}}^{\frac{6-r}{2r}}\left\Vert \nabla
_{h}f\right\Vert _{L^{2}}^{\frac{r-2}{r}}\left\Vert \nabla f\right\Vert
_{L^{2}}^{\frac{r-2}{2r}},  \label{eq9}
\end{eqnarray}%
for every $f\in H^{1}(\mathbb{R}^{3})$ and $r\in \lbrack 2,6]$.

Now we are ready to present the proof of Theorem \ref{th1}.

\begin{pf}
Since the initial data $\left( u_{0},\theta _{0}\right) \in H^{1}(\mathbb{R}%
^{3})$ with $\mathrm{div}\,u_{0}=0$ in $\mathbb{R}^{3}$, there exists a
unique local strong solution $(u,\theta )$ of the 3D Boussinesq equations on
$(0,T)$ (see \cite{AH, CN, CKN, HL}). By using a standard method, we only
need to to show the following a priori estimate
\begin{eqnarray}
&&\underset{0\leq t\leq T}{\sup }\left( \left\Vert \nabla u(\cdot
,t)\right\Vert _{L^{2}}^{2}+\left\Vert \nabla \theta (\cdot ,t)\right\Vert
_{L^{2}}^{2}\right)  \label{eq01} \\
&\leq &\left( \left\Vert \nabla u_{0}\right\Vert _{L^{2}}^{2}+\left\Vert
\nabla \theta _{0}\right\Vert _{L^{2}}^{2}+C\left\Vert \nabla
_{h}u_{0}\right\Vert _{L^{2}}^{\frac{8}{3}}+C\left\Vert \nabla _{h}\theta
_{0}\right\Vert _{L^{2}}^{\frac{8}{3}}\right) e^{C\mathcal{K}(T)},  \notag
\end{eqnarray}%
where
\begin{equation*}
\mathcal{K}(T)=\int\limits_{0}^{T}\left( \left\Vert \nabla
_{h}u(s)\right\Vert _{\overset{\cdot }{B}_{\infty ,\infty }^{-1}}^{\frac{8}{3%
}}+\left\Vert \nabla _{h}\theta (s)\right\Vert _{\overset{\cdot }{B}_{\infty
,\infty }^{-1}}^{\frac{8}{3}}\right) ds.
\end{equation*}%
Let
\begin{eqnarray*}
\mathcal{J}(t) &=&\left\Vert \nabla _{h}u(t)\right\Vert
_{L^{2}}^{2}+\left\Vert \nabla _{h}\theta (t)\right\Vert
_{L^{2}}^{2}+\int\limits_{0}^{t}\left( \left\Vert \nabla \nabla _{h}u(\tau
)\right\Vert _{L^{2}}^{2}+\left\Vert \nabla \nabla _{h}\theta (\tau
)\right\Vert _{L^{2}}^{2}\right) d\tau , \\
\mathcal{X}(t) &=&\int\limits_{0}^{t}\left( \left\Vert \nabla _{h}u(\tau
)\right\Vert _{\overset{\cdot }{B}_{\infty ,\infty }^{-1}}^{2}+\left\Vert
\nabla _{h}\theta (\tau )\right\Vert _{\overset{\cdot }{B}_{\infty ,\infty
}^{-1}}^{2}\right) \left( \left\Vert \nabla u(\tau )\right\Vert
_{L^{2}}^{2}+\left\Vert \nabla \theta (\tau )\right\Vert _{L^{2}}^{2}\right)
d\tau , \\
\mathcal{Z}(t) &=&\left\Vert \nabla u(\cdot ,t)\right\Vert
_{L^{2}}^{2}+\left\Vert \nabla \theta (\cdot ,t)\right\Vert
_{L^{2}}^{2}+\int\limits_{0}^{t}\left( \left\Vert \Delta u(\cdot ,\tau
)\right\Vert _{L^{2}}^{2}+\left\Vert \Delta \theta (\cdot ,\tau )\right\Vert
_{L^{2}}^{2}\right) d\tau , \\
\mathcal{W}(t) &=&\int\limits_{0}^{t}\left( \left\Vert \nabla _{h}u(\tau
)\right\Vert _{\overset{\cdot }{B}_{\infty ,\infty }^{-1}}^{\frac{8}{3}%
}+\left\Vert \nabla _{h}\theta (\tau )\right\Vert _{\overset{\cdot }{B}%
_{\infty ,\infty }^{-1}}^{\frac{8}{3}}\right) \left( \left\Vert \nabla
u(\tau )\right\Vert _{L^{2}}^{2}+\left\Vert \nabla \theta (\tau )\right\Vert
_{L^{2}}^{2}\right) d\tau .
\end{eqnarray*}

We start with the estimates of $\left\Vert \nabla _{h}u\right\Vert _{L^{2}}$
and $\left\Vert \nabla _{h}\theta \right\Vert _{L^{2}}$. Multiplying the
first equation of (\ref{eq1.1}) by $(-\Delta _{h}u)$, where $\Delta
_{h}=\partial _{1}\partial _{1}+\partial _{2}\partial _{2}$ is the
horizontal Laplacian and integrating by parts and using the divergence free
condition $\nabla \cdot u=0$ into account, we get%
\begin{eqnarray}
\frac{1}{2}\frac{d}{dt}\left\Vert \nabla _{h}u(t)\right\Vert
_{L^{2}}^{2}+\left\Vert \nabla \nabla _{h}u\right\Vert _{L^{2}}^{2} &=&\int_{%
\mathbb{R}^{3}}(u\cdot \nabla )u\cdot \Delta _{h}udx-\int_{\mathbb{R}%
^{3}}\theta e_{3}\cdot \Delta _{h}udx  \notag \\
&=&\int_{\mathbb{R}^{3}}\sum_{j=1}^{3}\sum_{l=1}^{2}u_{j}\partial _{j}u\cdot
\partial _{l}^{2}udx-\int_{\mathbb{R}^{3}}\sum_{l=1}^{2}\theta e_{3}\cdot
\partial _{l}^{2}udx  \notag \\
&=&-\int_{\mathbb{R}^{3}}\sum_{j=1}^{3}\sum_{l=1}^{2}\partial
_{l}u_{j}\partial _{j}u\partial _{l}udx+\int_{\mathbb{R}^{3}}\sum_{l=1}^{2}%
\partial _{l}(\theta e_{3})\partial _{l}udx  \notag \\
&=&-\int_{\mathbb{R}^{3}}\nabla _{h}u\cdot \nabla u\cdot \nabla
_{h}udx+\int_{\mathbb{R}^{3}}\nabla _{h}(\theta e_{3})\cdot \nabla _{h}udx,
\label{eq21}
\end{eqnarray}%
where we have used%
\begin{equation*}
\int_{\mathbb{R}^{3}}\sum_{j=1}^{3}\sum_{l=1}^{2}u_{j}\partial _{j}\partial
_{l}u\cdot \partial _{l}udx=0.
\end{equation*}
Similarly, multiplying the second equation of (\ref{eq1.1}) by $(-\Delta
_{h}\theta )$, we obtain%
\begin{eqnarray}
\frac{1}{2}\frac{d}{dt}\left\Vert \nabla _{h}\theta (t)\right\Vert
_{L^{2}}^{2}+\left\Vert \nabla \nabla _{h}\theta \right\Vert _{L^{2}}^{2}
&=&\int_{\mathbb{R}^{3}}(u\cdot \nabla )\theta \cdot \Delta _{h}\theta
dx=\int_{\mathbb{R}^{3}}\sum_{j=1}^{3}\sum_{l=1}^{2}u_{j}\partial _{j}\theta
\cdot \partial _{l}^{2}\theta dx  \notag \\
&=&-\int_{\mathbb{R}^{3}}\sum_{j=1}^{3}\sum_{l=1}^{2}\partial
_{l}u_{j}\partial _{j}\theta \partial _{l}\theta dx=-\int_{\mathbb{R}%
^{3}}\nabla _{h}u\cdot \nabla \theta \cdot \nabla _{h}\theta dx,
\label{eq22}
\end{eqnarray}%
wher we have used%
\begin{equation*}
\int_{\mathbb{R}^{3}}\sum_{j=1}^{3}\sum_{l=1}^{2}u_{j}\partial _{j}\partial
_{l}\theta \cdot \partial _{l}\theta dx=0.
\end{equation*}
Combining (\ref{eq21}) and (\ref{eq22}) yields
\begin{eqnarray}
\lefteqn{\frac{1}{2}\frac{d}{dt}(\left\Vert \nabla _{h}u(t)\right\Vert
_{L^{2}}^{2}+\left\Vert \nabla _{h}\theta (t)\right\Vert
_{L^{2}}^{2})+\left\Vert \nabla \nabla _{h}u\right\Vert
_{L^{2}}^{2}+\left\Vert \nabla \nabla _{h}\theta \right\Vert _{L^{2}}^{2}}
\notag \\
&=&-\int_{\mathbb{R}^{3}}\nabla _{h}u\cdot \nabla u\cdot \nabla
_{h}udx-\int_{\mathbb{R}^{3}}\nabla _{h}u\cdot \nabla \theta \cdot \nabla
_{h}\theta dx+\int_{\mathbb{R}^{3}}\nabla _{h}(\theta e_{3})\cdot \nabla
_{h}udx  \notag \\
&=&R_{1}+R_{2}+R_{3},  \label{eq3.4}
\end{eqnarray}%
Atention is now focused on bounding these terms; we start with $R_{1}$.
Using H\"{o}lder and Young's inequalities, one has, for $R_{1},$ \
\begin{eqnarray*}
\left\vert R_{1}\right\vert &\leq &C\left\Vert \nabla _{h}u\right\Vert
_{L^{4}}^{2}\left\Vert \nabla u\right\Vert _{L^{2}} \\
&\leq &C\left\Vert \nabla \nabla _{h}u\right\Vert _{L^{2}}\left\Vert \nabla
_{h}u\right\Vert _{\overset{\cdot }{B}_{\infty ,\infty }^{-1}}\left\Vert
\nabla u\right\Vert _{L^{2}} \\
&\leq &\frac{1}{4}\left\Vert \nabla \nabla _{h}u\right\Vert
_{L^{2}}^{2}+C\left\Vert \nabla _{h}u\right\Vert _{\overset{\cdot }{B}%
_{\infty ,\infty }^{-1}}^{2}\left\Vert \nabla u\right\Vert _{L^{2}}^{2}.
\end{eqnarray*}%
Here we have used the following inequality due to Meyer-Gerard-Oru \cite{MGO}
(see also \cite{GG}) :
\begin{equation}
\left\Vert f\right\Vert _{L^{4}}^{2}\leq C\left\Vert \nabla f\right\Vert
_{L^{2}}\left\Vert f\right\Vert _{\overset{\cdot }{B}_{\infty ,\infty
}^{-1}}.  \label{eq69}
\end{equation}%
For $R_{2}$, analogously, using H\"{o}lder and Young's inequalities, we
deduce from (\ref{eq69}) that \
\begin{eqnarray*}
\left\vert R_{2}\right\vert &\leq &C\left\Vert \nabla _{h}u\right\Vert
_{L^{4}}\left\Vert \nabla \theta \right\Vert _{L^{2}}\left\Vert \nabla
_{h}\theta \right\Vert _{L^{4}} \\
&\leq &C\left\Vert \nabla \nabla _{h}u\right\Vert _{L^{2}}^{\frac{1}{2}%
}\left\Vert \nabla _{h}u\right\Vert _{\overset{\cdot }{B}_{\infty ,\infty
}^{-1}}^{\frac{1}{2}}\left\Vert \nabla \theta \right\Vert _{L^{2}}\left\Vert
\nabla \nabla _{h}\theta \right\Vert _{L^{2}}^{\frac{1}{2}}\left\Vert \nabla
_{h}\theta \right\Vert _{\overset{\cdot }{B}_{\infty ,\infty }^{-1}}^{\frac{1%
}{2}} \\
&=&C\left( \left\Vert \nabla \nabla _{h}u\right\Vert _{L^{2}}^{2}\right) ^{%
\frac{1}{4}}\left( \left\Vert \nabla _{h}u\right\Vert _{\overset{\cdot }{B}%
_{\infty ,\infty }^{-1}}^{2}\left\Vert \nabla \theta \right\Vert
_{L^{2}}^{2}\right) ^{\frac{1}{4}}\left( \left\Vert \nabla \nabla _{h}\theta
\right\Vert _{L^{2}}^{2}\right) ^{\frac{1}{4}}\left( \left\Vert \nabla
_{h}\theta \right\Vert _{\overset{\cdot }{B}_{\infty ,\infty
}^{-1}}^{2}\left\Vert \nabla \theta \right\Vert _{L^{2}}^{2}\right) ^{\frac{1%
}{4}} \\
&\leq &\frac{1}{4}\left\Vert \nabla \nabla _{h}u\right\Vert _{L^{2}}^{2}+%
\frac{1}{4}\left\Vert \nabla \nabla _{h}\theta \right\Vert
_{L^{2}}^{2}+C\left\Vert \nabla _{h}u\right\Vert _{\overset{\cdot }{B}%
_{\infty ,\infty }^{-1}}^{2}\left\Vert \nabla \theta \right\Vert
_{L^{2}}^{2}+C\left\Vert \nabla _{h}\theta \right\Vert _{\overset{\cdot }{B}%
_{\infty ,\infty }^{-1}}^{2}\left\Vert \nabla \theta \right\Vert
_{L^{2}}^{2}.
\end{eqnarray*}%
For $R_{3}$, by means of the H\"{o}lder and Cauchy inequalities, it follows
that%
\begin{eqnarray*}
\left\vert R_{3}\right\vert &\leq &C\left\Vert \nabla _{h}\theta \right\Vert
_{L^{2}}\left\Vert \nabla _{h}u\right\Vert _{L^{2}} \\
&\leq &C(\left\Vert \nabla _{h}\theta \right\Vert _{L^{2}}^{2}+\left\Vert
\nabla _{h}u\right\Vert _{L^{2}}^{2}).
\end{eqnarray*}%
Inserting the above estimate into (\ref{eq3.4}), we derive that
\begin{eqnarray*}
\lefteqn{\frac{d}{dt}\left( \left\Vert \nabla _{h}u(\cdot ,t)\right\Vert
_{L^{2}}^{2}+\left\Vert \nabla _{h}\theta (\cdot ,t)\right\Vert
_{L^{2}}^{2}\right) +\left\Vert \nabla \nabla _{h}u(\cdot ,t)\right\Vert
_{L^{2}}^{2}+\left\Vert \nabla \nabla _{h}\theta (\cdot ,t)\right\Vert
_{L^{2}}^{2}} \\
&\leq &C\left( \left\Vert \nabla _{h}u(\cdot ,t)\right\Vert _{\overset{\cdot
}{B}_{\infty ,\infty }^{-1}}^{2}+\left\Vert \nabla _{h}\theta (\cdot
,t)\right\Vert _{\overset{\cdot }{B}_{\infty ,\infty }^{-1}}^{2}\right)
\left( \left\Vert \nabla u(\cdot ,t)\right\Vert _{L^{2}}^{2}+\left\Vert
\nabla \theta (\cdot ,t)\right\Vert _{L^{2}}^{2}\right) .
\end{eqnarray*}%
Integrating the above inequality in time variable over $0\leq \tau \leq t$,
we get
\begin{equation*}
\mathcal{J}(t)\leq \left\Vert \nabla _{h}u_{0}\right\Vert
_{L^{2}}^{2}+\left\Vert \nabla _{h}\theta _{0}\right\Vert _{L^{2}}^{2}+C%
\mathcal{X}(t).
\end{equation*}

Next, we derive the bounds of $\left\Vert \nabla u\right\Vert _{L^{2}}$ and $%
\left\Vert \nabla \theta \right\Vert _{L^{2}}$. Multiplying the two
equations of (\ref{eq1.1}) by $(-\Delta u)$ and $(-\Delta \theta )$,
repectively, integrating and applying the incompressibility condition, we
have by (\ref{eq20.1})
\begin{eqnarray}
\lefteqn{\frac{1}{2}\frac{d}{dt}(\left\Vert \nabla u(t)\right\Vert
_{L^{2}}^{2}+\left\Vert \nabla \theta (t)\right\Vert
_{L^{2}}^{2})+\left\Vert \Delta u(t)\right\Vert _{L^{2}}^{2}+\left\Vert
\Delta \theta (t)\right\Vert _{L^{2}}^{2}}  \notag \\
&=&\int_{\mathbb{R}^{3}}(u\cdot \nabla )u\cdot \Delta udx+\int_{\mathbb{R}%
^{3}}(u\cdot \nabla )\theta \cdot \Delta \theta dx+\int_{\mathbb{R}%
^{3}}(\theta e_{3})\cdot \Delta udx  \notag \\
&\leq &C\int_{\mathbb{R}^{3}}\left\vert \nabla _{h}u\right\vert \left\vert
\nabla u\right\vert ^{2}dx+C\int_{\mathbb{R}^{3}}\left\vert \nabla
_{h}u\right\vert \left\vert \nabla \theta \right\vert ^{2}dx+C\int_{\mathbb{R%
}^{3}}\left\vert \nabla _{h}\theta \right\vert \left\vert \nabla
u\right\vert \left\vert \nabla \theta \right\vert dx+C\left\Vert \theta
\right\Vert _{L^{2}}\left\Vert \Delta u\right\Vert _{L^{2}}  \notag \\
&\leq &K_{1}+K_{2}+K_{3}+\frac{1}{4}\left\Vert \Delta u\right\Vert
_{L^{2}}^{2}+C\left\Vert \theta \right\Vert _{L^{2}}^{2}.  \label{eq0}
\end{eqnarray}%
Now we deal with $K_{1}$. It \ follows that, from the H\"{o}lder inequality
and (\ref{eq9})
\begin{eqnarray*}
K_{1} &\leq &C\left\Vert \nabla _{h}u\right\Vert _{L^{2}}\left\Vert \nabla
u\right\Vert _{L^{4}}^{2} \\
&\leq &C\left\Vert \nabla _{h}u\right\Vert _{L^{2}}\left\Vert \nabla
u\right\Vert _{L^{2}}^{\frac{1}{2}}\left\Vert \nabla \nabla _{h}u\right\Vert
_{L^{2}}\left\Vert \Delta u\right\Vert _{L^{2}}^{\frac{1}{2}}.
\end{eqnarray*}%
For $K_{2}$, H\"{o}lder inequality and (\ref{eq9}), together give,
\begin{eqnarray*}
K_{2} &\leq &C\left\Vert \nabla _{h}u\right\Vert _{L^{2}}\left\Vert \nabla
\theta \right\Vert _{L^{4}}^{2} \\
&\leq &C\left\Vert \nabla _{h}u\right\Vert _{L^{2}}\left\Vert \nabla \theta
\right\Vert _{L^{2}}^{\frac{1}{2}}\left\Vert \nabla \nabla _{h}\theta
\right\Vert _{L^{2}}\left\Vert \Delta \theta \right\Vert _{L^{2}}^{\frac{1}{2%
}}.
\end{eqnarray*}%
Arguing similarly as the estimate of $K_{1}$, thanks to the H\"{o}lder
inequality and (\ref{eq9}), one has%
\begin{eqnarray*}
K_{3} &\leq &C\left\Vert \nabla _{h}\theta \right\Vert _{L^{2}}\left\Vert
\nabla u\right\Vert _{L^{4}}\left\Vert \nabla \theta \right\Vert _{L^{4}} \\
&\leq &C\left\Vert \nabla _{h}\theta \right\Vert _{L^{2}}\left( \left\Vert
\nabla u\right\Vert _{L^{2}}^{\frac{1}{4}}\left\Vert \nabla \nabla
_{h}u\right\Vert _{L^{2}}^{\frac{1}{2}}\left\Vert \Delta u\right\Vert
_{L^{2}}^{\frac{1}{4}}\right) \left( \left\Vert \nabla \theta \right\Vert
_{L^{2}}^{\frac{1}{4}}\left\Vert \nabla \nabla _{h}\theta \right\Vert
_{L^{2}}^{\frac{1}{2}}\left\Vert \Delta \theta \right\Vert _{L^{2}}^{\frac{1%
}{4}}\right) .
\end{eqnarray*}%
Combining the above estimates of $K_{1},K_{2}$ and $K_{3}$ and inserting
into (\ref{eq0}), we get
\begin{eqnarray*}
\lefteqn{\frac{d}{dt}(\left\Vert \nabla u(\cdot ,t)\right\Vert
_{L^{2}}^{2}+\left\Vert \nabla \theta (\cdot ,t)\right\Vert
_{L^{2}}^{2})+\left\Vert \Delta u\right\Vert _{L^{2}}^{2}+\left\Vert \Delta
\theta \right\Vert _{L^{2}}^{2}} \\
&\leq &C\left\Vert \nabla _{h}u\right\Vert _{L^{2}}\left\Vert \nabla
u\right\Vert _{L^{2}}^{\frac{1}{2}}\left\Vert \nabla \nabla _{h}u\right\Vert
_{L^{2}}\left\Vert \Delta u\right\Vert _{L^{2}}^{\frac{1}{2}}+C\left\Vert
\nabla _{h}u\right\Vert _{L^{2}}\left\Vert \nabla \theta \right\Vert
_{L^{2}}^{\frac{1}{2}}\left\Vert \nabla \nabla _{h}\theta \right\Vert
_{L^{2}}\left\Vert \Delta \theta \right\Vert _{L^{2}}^{\frac{1}{2}} \\
&&+C\left\Vert \nabla _{h}\theta \right\Vert _{L^{2}}\left\Vert \nabla
u\right\Vert _{L^{2}}^{\frac{1}{4}}\left\Vert \nabla \nabla _{h}u\right\Vert
_{L^{2}}^{\frac{1}{2}}\left\Vert \Delta u\right\Vert _{L^{2}}^{\frac{1}{4}%
}\left\Vert \nabla \theta \right\Vert _{L^{2}}^{\frac{1}{4}}\left\Vert
\nabla \nabla _{h}\theta \right\Vert _{L^{2}}^{\frac{1}{2}}\left\Vert \Delta
\theta \right\Vert _{L^{2}}^{\frac{1}{4}}+C\left\Vert \theta \right\Vert
_{L^{2}}^{2}.
\end{eqnarray*}%
Integrating the above inequality in time variable over $0\leq \tau \leq t$,
one shows that
\begin{eqnarray*}
\mathcal{Z}(t) &\leq &C\left( \underset{0\leq \tau \leq t}{\sup }\left\Vert
\nabla _{h}u\right\Vert _{L^{2}}\right) A_{1}(t)+C\left( \underset{0\leq
\tau \leq t}{\sup }\left\Vert \nabla _{h}u\right\Vert _{L^{2}}\right)
A_{2}(t) \\
&&+C\left( \underset{0\leq \tau \leq t}{\sup }\left\Vert \nabla _{h}\theta
\right\Vert _{L^{2}}\right) \left( A_{3}(t)\times A_{4}(t)\right)
+\left\Vert \nabla u_{0}\right\Vert _{L^{2}}^{2}+\left\Vert \nabla \theta
_{0}\right\Vert _{L^{2}}^{2} \\
&\leq &C\left( \left\Vert \nabla _{h}u_{0}\right\Vert
_{L^{2}}^{2}+\left\Vert \nabla _{h}\theta _{0}\right\Vert _{L^{2}}^{2}+%
\mathcal{X}(t)\right) \times A_{5}(t)+\left\Vert \nabla u_{0}\right\Vert
_{L^{2}}^{2}+\left\Vert \nabla \theta _{0}\right\Vert _{L^{2}}^{2},
\end{eqnarray*}%
where we denote%
\begin{eqnarray*}
A_{1}(t) &=&\left( \int\limits_{0}^{t}\left\Vert \nabla u(\tau )\right\Vert
_{L^{2}}^{2}d\tau \right) ^{\frac{1}{4}}\left(
\int\limits_{0}^{t}\left\Vert \nabla \nabla _{h}u(\tau )\right\Vert
_{L^{2}}^{2}d\tau \right) ^{\frac{1}{2}}\left(
\int\limits_{0}^{t}\left\Vert \Delta u(\tau )\right\Vert _{L^{2}}^{2}d\tau
\right) ^{\frac{1}{4}}, \\
A_{2}(t) &=&\left( \int\limits_{0}^{t}\left\Vert \nabla \theta (\tau
)\right\Vert _{L^{2}}^{2}d\tau \right) ^{\frac{1}{4}}\left(
\int\limits_{0}^{t}\left\Vert \nabla \nabla _{h}\theta (\tau )\right\Vert
_{L^{2}}^{2}d\tau \right) ^{\frac{1}{2}}\left(
\int\limits_{0}^{t}\left\Vert \Delta \theta (\tau )\right\Vert
_{L^{2}}^{2}d\tau \right) ^{\frac{1}{4}}, \\
A_{3}(t) &=&\left( \int\limits_{0}^{t}\left\Vert \nabla u(\tau )\right\Vert
_{L^{2}}^{2}d\tau \right) ^{\frac{1}{8}}\left(
\int\limits_{0}^{t}\left\Vert \nabla \nabla _{h}u(\tau )\right\Vert
_{L^{2}}^{2}d\tau \right) ^{\frac{1}{4}}\left(
\int\limits_{0}^{t}\left\Vert \Delta u(\tau )\right\Vert _{L^{2}}^{2}d\tau
\right) ^{\frac{1}{8}}, \\
A_{4}(t) &=&\left( \int\limits_{0}^{t}\left\Vert \nabla \theta (\tau
)\right\Vert _{L^{2}}^{2}d\tau \right) ^{\frac{1}{8}}\left(
\int\limits_{0}^{t}\left\Vert \nabla \nabla _{h}\theta (\tau )\right\Vert
_{L^{2}}^{2}d\tau \right) ^{\frac{1}{4}}\left(
\int\limits_{0}^{t}\left\Vert \Delta \theta (\tau )\right\Vert
_{L^{2}}^{2}d\tau \right) ^{\frac{1}{8}}, \\
A_{5}(t) &=&\left( \int\limits_{0}^{t}\left( \left\Vert \Delta u(\cdot
,\tau )\right\Vert _{L^{2}}^{2}+\left\Vert \Delta \theta (\cdot ,\tau
)\right\Vert _{L^{2}}^{2}\right) d\tau \right) ^{\frac{1}{4}}.
\end{eqnarray*}

By virtue of the H\"{o}lder and Young inequalities and energy inequality (%
\ref{eq100}), we get
\begin{eqnarray*}
\mathcal{Z}(t) &\leq &\left\Vert \nabla u_{0}\right\Vert
_{L^{2}}^{2}+\left\Vert \nabla \theta _{0}\right\Vert _{L^{2}}^{2}+C\left(
\left\Vert \nabla _{h}u_{0}\right\Vert _{L^{2}}^{\frac{8}{3}}+\left\Vert
\nabla _{h}\theta _{0}\right\Vert _{L^{2}}^{\frac{8}{3}}\right) \\
&&+C\int\limits_{0}^{t}\left( \left\Vert \nabla _{h}u(\tau )\right\Vert _{%
\overset{\cdot }{B}_{\infty ,\infty }^{-1}}^{2}\left\Vert \nabla u(\tau
)\right\Vert _{L^{2}}^{2}d\tau \right) ^{\frac{4}{3}}+C\int\limits_{0}^{t}%
\left( \left\Vert \nabla _{h}u(\tau )\right\Vert _{\overset{\cdot }{B}%
_{\infty ,\infty }^{-1}}^{2}\left\Vert \nabla \theta (\tau )\right\Vert
_{L^{2}}^{2}d\tau \right) ^{\frac{4}{3}} \\
&&+C\int\limits_{0}^{t}\left( \left\Vert \nabla _{h}\theta (\tau
)\right\Vert _{\overset{\cdot }{B}_{\infty ,\infty }^{-1}}^{2}\left\Vert
\nabla \theta (\tau )\right\Vert _{L^{2}}^{2}d\tau \right) ^{\frac{4}{3}%
}+C\int\limits_{0}^{t}\left( \left\Vert \nabla _{h}\theta (\tau
)\right\Vert _{\overset{\cdot }{B}_{\infty ,\infty }^{-1}}^{2}\left\Vert
\nabla u(\tau )\right\Vert _{L^{2}}^{2}d\tau \right) ^{\frac{4}{3}} \\
&\leq &\left\Vert \nabla u_{0}\right\Vert _{L^{2}}^{2}+\left\Vert \nabla
\theta _{0}\right\Vert _{L^{2}}^{2}+C\left\Vert \nabla _{h}u_{0}\right\Vert
_{L^{2}}^{\frac{8}{3}}+C\left\Vert \nabla _{h}\theta _{0}\right\Vert
_{L^{2}}^{\frac{8}{3}} \\
&&+C\int\limits_{0}^{t}\left( \left\Vert \nabla _{h}u(\tau )\right\Vert _{%
\overset{\cdot }{B}_{\infty ,\infty }^{-1}}^{\frac{8}{3}}\left\Vert \nabla
u(\tau )\right\Vert _{L^{2}}^{2}d\tau \right) \left(
\int\limits_{0}^{t}\left\Vert \nabla u(\tau )\right\Vert _{L^{2}}^{2}d\tau
\right) ^{\frac{1}{3}} \\
&&+C\int\limits_{0}^{t}\left( \left\Vert \nabla _{h}u(\tau )\right\Vert _{%
\overset{\cdot }{B}_{\infty ,\infty }^{-1}}^{\frac{8}{3}}\left\Vert \nabla
\theta (\tau )\right\Vert _{L^{2}}^{2}d\tau \right) \left(
\int\limits_{0}^{t}\left\Vert \nabla \theta (\tau )\right\Vert
_{L^{2}}^{2}d\tau \right) ^{\frac{1}{3}} \\
&&+C\int\limits_{0}^{t}\left( \left\Vert \nabla _{h}\theta (\tau
)\right\Vert _{\overset{\cdot }{B}_{\infty ,\infty }^{-1}}^{\frac{8}{3}%
}\left\Vert \nabla \theta (\tau )\right\Vert _{L^{2}}^{2}d\tau \right)
\left( \int\limits_{0}^{t}\left\Vert \nabla \theta (\tau )\right\Vert
_{L^{2}}^{2}d\tau \right) ^{\frac{1}{3}} \\
&&+C\int\limits_{0}^{t}\left( \left\Vert \nabla _{h}\theta (\tau
)\right\Vert _{\overset{\cdot }{B}_{\infty ,\infty }^{-1}}^{\frac{8}{3}%
}\left\Vert \nabla u(\tau )\right\Vert _{L^{2}}^{2}d\tau \right) \left(
\int\limits_{0}^{t}\left\Vert \nabla u(\tau )\right\Vert _{L^{2}}^{2}d\tau
\right) ^{\frac{1}{3}}.
\end{eqnarray*}%
Thanks to the energy inequality (\ref{eq100}), we get
\begin{equation*}
\mathcal{Z}(t)\leq \left\Vert \nabla u_{0}\right\Vert
_{L^{2}}^{2}+\left\Vert \nabla \theta _{0}\right\Vert
_{L^{2}}^{2}+C\left\Vert \nabla _{h}u_{0}\right\Vert _{L^{2}}^{\frac{8}{3}%
}+C\left\Vert \nabla _{h}\theta _{0}\right\Vert _{L^{2}}^{\frac{8}{3}}+C%
\mathcal{W}(t).
\end{equation*}%
Taking the Gronwall inequality into consideration, we arrive at
\begin{equation*}
\mathcal{Z}(t)\leq \left( \left\Vert \nabla u_{0}\right\Vert
_{L^{2}}^{2}+\left\Vert \nabla \theta _{0}\right\Vert
_{L^{2}}^{2}+C\left\Vert \nabla _{h}u_{0}\right\Vert _{L^{2}}^{\frac{8}{3}%
}+C\left\Vert \nabla _{h}\theta _{0}\right\Vert _{L^{2}}^{\frac{8}{3}%
}\right) e^{C\mathcal{K}(t)},
\end{equation*}%
for any $t\in \lbrack 0,T)$, which implies that%
\begin{equation}
\underset{0\leq t\leq T}{\sup }\left( \left\Vert \nabla u(\cdot
,t)\right\Vert _{L^{2}}^{2}+\left\Vert \nabla \theta (\cdot ,t)\right\Vert
_{L^{2}}^{2}\right) <\infty .  \label{eq19}
\end{equation}%
In the end, by the standard arguments of continuation of local solutions, we
complete the proof of Theorem \ref{th1}.
\end{pf}

%

\end{document}